\documentclass[1p,preprint,12pt]{elsarticle}

\usepackage{graphicx}
\usepackage{epstopdf}
\usepackage{subcaption}
\usepackage{algorithm}
\usepackage[noend]{algpseudocode}
\usepackage{amssymb}
\usepackage{bm}
\usepackage{amsmath,amsthm}
\usepackage{xmpmulti}
\usepackage{multicol}
\usepackage{mathtools}
\usepackage{tabularx}
\usepackage{color,soul,colortbl} 
\usepackage{tikz}
\usetikzlibrary{fit}
\usetikzlibrary{shapes}
\usetikzlibrary{arrows}
\usetikzlibrary{decorations.markings}
\usetikzlibrary{positioning}

\tikzset{mylabel/.style={font=\footnotesize}}
\tikzset{mymidlabel/.style={fill=white,font=\footnotesize}}
\newcommand{\bsub}{\begin{subequations}}

\newcommand{\mbr}{\mathbb{R}}

\def \vx{\mbox{\boldmath $x$}}
\def \vn{\mbox{\boldmath $n$}}
\def \vxi{\mbox{\boldmath $\xi$}}

\def \nind{\noindent}

\def \vs{\vspace{.2in}}

\newcommand{\ds}{\displaystyle}
\newcommand{\efig}{\end{figure}}
\newcommand{\bfig}{\begin{figure}}
\newcommand{\bea}{\begin{eqnarray}}
\newcommand{\eea}{\end{eqnarray}}

\newcommand{\beq}{\begin{equation}}
\newcommand{\eeq}{\end{equation}}
\newcommand{\ba}{\begin{array}}
\newcommand{\ea}{\end{array}}

\newcommand{\fr}{\frac}

\usepackage{diagbox,multirow}
\newcolumntype{L}[1]{>{\raggedright\let\newline\\\arraybackslash\hspace{0pt}}m{#1}}
\newcolumntype{C}[1]{>{\centering\arraybackslash}p{#1}}	% center the text in cells of table with fixed width
\usepackage[T1]{fontenc}
\usepackage[font=small,labelfont=bf,tableposition=top]{caption}

\DeclareCaptionLabelFormat{andtable}{#1~#2  \&  \tablename~\thetable}
\usepackage{threeparttable}

\usepackage[pdftex]{hyperref}

\journal{J. Comput. Math.  Appl.}
\begin{document}

\begin{frontmatter}



\title{Domain Decomposition Parabolic Monge-Amp\`ere Approach for Fast Generation of Adaptive Moving Meshes} 


\author[aa]{M.~Sulman\corref{cor1}}
\ead{mohamed.sulman@wright.edu}
\cortext[cor1]{Corresponding author}

\author[bb]{T.~Nguyen}
\ead{tbnguyen@lanl.gov}

\author[cc]{R.~Haynes}
\ead{rhaynes@mun.ca}

\author[dd]{W.~Huang}
\ead{whuang@ku.edu}

\address[aa]{\small Department of Mathematics \& Statistics, Wright State University, Dayton, OH 45435, USA}

\address[bb]{\small Los Alamos National Laboratory, Los Alamos, NM 87545, USA}

\address[cc]{\small Department of Mathematics \& Statistics, Memorial University of Newfoundland, NL A1C 5S7, CANADA}

\address[dd]{\small Department of Mathematics, University of Kansas, Lawrence, KS 66045, USA} 

\begin{abstract}

A fast method is presented for adaptive moving mesh generation in multi-dimensions using a domain decomposition parabolic Monge-Amp\`ere approach.  The domain decomposition procedure employed here  is non-iterative and involves splitting the computational domain into overlapping subdomains. An adaptive mesh on each subdomain is then computed as the image of the solution of the $L^2$ optimal mass transfer problem using a parabolic Monge-Amp\`ere method. The  domain decomposition approach allows straightforward implementation for the parallel computation of adaptive meshes which helps to reduce computational time significantly.  Results are presented to show   
the numerical convergence of the domain decomposition solution to the single domain solution. Several numerical experiments are given to demonstrate  the performance and  efficiency of the proposed  method.
The numerical results indicate that the domain decomposition parabolic Monge-Amp\`ere method is more efficient than the standard implementation of the parabolic Monge-Amp\`ere method on the whole domain, in particular when computing adaptive meshes in three spatial dimensions. 

   \end{abstract}

\end{frontmatter}



\section{Introduction} \label{Sect1}

Adaptive mesh methods have become  increasingly popular over the last three decades.    
The use of a uniform grid for solving partial differential equations (PDEs) can be prohibitively expensive, especially for  problems  in multi-dimensions  where their solutions develop sharp structures in some small regions of the physical domain. The large errors in the approximation of the physical solution are expected to occur in those regions. Therefore, mesh adaptation is needed to improve  the accuracy of the  numerical solution  while reducing the computational cost; e.g. see \cite{Dorfi:1987, Thompson:1985, Anderson:1987, Huang2011}.    Numerous adaptive mesh techniques have been developed over the last three decades  (see, for example,  \cite{Miller:1981, Millera:1981, Gelinas:1981, Adjerid:1986, Furrzeland:1990, Hawken:1991, Adjerid:1992, Huang:1994, Huang:1999, Huang:2003, Cao:2002, Backett:2002}).  Here, we consider a class of adaptive mesh methods called the $r$-refinement or moving mesh method.  In this approach the mesh points are continuously redistributed so that they are concentrated in the  regions of large solution variations or gradients. In one spatial dimension, the adaptive mesh can be computed based on the equidistribution principle \cite{Boor:1973, Huang:1994b}; an estimate of the numerical approximation error is evenly distributed among 
mesh elements. The equidistribution condition alone is insufficient to uniquely determine an adaptive mesh in multi-dimensions. A number of adaptive moving mesh methods have been developed so far; for example,
see \cite{Huang2011,Tan05,BHR09} and references therein.
In this work we are interested in methods based on solving the optimal mass transfer problem \cite{Budd:2009, Sulman:2011a, Chacon:2011}. The optimal mass transfer problem, also known as the  Monge-Kantorovich problem (MKP), appears in  numerous applications in science and engineering  \cite{Monge:1781, Kantorovich:1948, Benamou:2000}. Sulman et al.  \cite{Sulman:2011a, Sulman:2011b}  describe finding the optimal solution of the $L^2$ MKP as the steady state solution of a parabolic Monge-Amp\`ere equation (PMA). 

In this paper, we present a domain decomposition parabolic Monge-Amp\`ere 
(DDPMA)  moving  mesh method for the generation of adaptive meshes in multi-dimensions. There are two main advantages of domain decomposition approaches for solving PDEs. First, DD is a natural approach for computing the  numerical solutions of both steady-state and time-dependent PDEs in parallel. Second, DD allows the use of different time steps on different subdomains in the time dependent context. These ideas can significantly reduce the computational time for adaptive mesh generation.  
In one dimension DD has been studied theoretically for mesh generation based on the equidistribution principle at the continuous level in \cite{Haynes:2012} and at the discrete level in \cite{Haynes:2017}.  Numerical DD methods for PDE based mesh generation in multi-dimensions can also be found in \cite{Haynes:2014, Bihlo:2014}.
A recent summary is available in \cite{Haynes:2018}.
Our approach here is based on solving the time dependent PMA equation. 
There are three common approaches used to apply domain decomposition to parabolic problems. The first approach is to apply the traditional iterative Schwarz algorithms first developed for elliptic problems \cite{Schwarz:1869, Lions:1989, Babuska:1958, Michlin:1951, White:1987, Jeltsch:1995} to the elliptic equations which arise upon semi-discretizing the time-dependent PDE in time (see \cite{Gander:1997, Gander:2002, Gander:2005, Vabishchevich:2008, Cai:1994, Zheng:2008, Qin:2008}).  The second approach is to split the whole space-time domain into overlapping or non-overlapping space-time subdomains in a Schwarz waveform relaxation framework \cite{Lions:1989, Mota:2017, Gnatyuk:2015, Dai:2016}.  The third approach is non-iterative domain decomposition which  is used to further reduce the computational cost  \cite{Dawsonfinitedifferencedomain1991, DawsonExplicitImplicitConservative1994, ZhangStablegloballynoniterative2000, GuangweiUNCONDITIONALSTABILITYPARALLEL2007, YangNoniterativeparallelSchwarz2017, Xuenewparallelalgorithm2018}.  Motivated by  this literature, the DDPMA 
  method proceeds by splitting the computational domain $\Omega_c$  block-wise or slab-wise into overlapping subdomains  and computes an approximation to the solution of the MA equation using a  non-iterative DD approach at each time level of  the pseudo time integration of the nonlinear parabolic Monge-Amp\`ere equation. We will study if this can be done without sacrificing mesh quality.

The paper is organized as follows. In Section~\ref{Sect2}, we give a brief description of the parabolic Monge-Amp\`ere method for generating adaptive meshes in multi-dimensions based on solving the $L^2$ optimal mass transfer problem. In Section~\ref{Sect3}, we describe the domain decomposition parabolic Monge-Amp\`ere moving  mesh method. In Section \ref{Sect4}, several numerical experiments  are presented to demonstrate the performance and efficiency of the proposed DDPMA method including results on the numerical convergence of the method.  Lastly, a discussion of  the results and  some concluding remarks are given  in Section \ref{Sect5}.

\section{The parabolic Monge-Amp\`ere adaptive mesh method} \label{Sect2}

 The parabolic Monge-Amp\`ere (PMA) method computes an adaptive mesh at any time $t$ as the image of a coordinate transformation   $ \vx = \vx(\vxi),$ defined from the logical or computational domain $\Omega_c\subset \mbr^d$ ($d \ge 1$) to the physical domain $\Omega\subset \mbr^d$.  The transformation $\vx = \vx(\vxi)$ is determined by equidistributing  a measure $\rho(\vx)$ of the solution error or variation  over mesh elements in the physical domain $\Omega$
\cite{Thompson:1985, Huang2011, Boor:1973, Sulman:2011a}. The equidistribution of  $\rho(\vx)$ can  be expressed \cite{Huang:1994, Boor:1973} by the constraint
 \beq\label{edp1}
\rho(\vx(\vxi)) \mbox{J}(\vx(\vxi))  = 1,  \quad \vxi \in \Omega_c, \;  \vx \in \Omega
\eeq

\nind where $\mbox{J}$ is the Jacobian of the transformation.  In two spatial dimensions, we have $\vx = (x,y)$, $\vxi = (\xi,\eta)$, and $\mbox{J} = x_{\xi} y_{\eta} - x_{\eta} y_{\xi}$.    

Notice that here, we require  that the mesh density function $\rho(\vx)$ to be normalized, i.e. $\int_\Omega \rho(\vx) d\vx = 1,$  and in this case the right hand side of (\ref{edp1}) will be modified to $\ds{1/\left|\Omega_c\right|.}$  Thus, the constraint (\ref{edp1}) takes the form 
\beq\label{edp}
\left|\Omega_c\right|\rho(\vx) \mbox{J}  = 1,  \quad \vxi \in \Omega_c, \;  \vx \in \Omega .
\eeq

The equidistribution constraint (\ref{edp}) alone is insufficient to uniquely determine the coordinate transformation, $\vx = \vx(\vxi )$  in multi-dimensions.  If the solution of the physical model does not involve large variations in the physical domain, then the spatial derivatives of physical solutions can be accurately approximated using standard finite difference schemes on a uniform grid. In this case the coordinate transformation    $\vx = \vx(\vxi )$ corresponds to the identity map, and  the constraint (\ref{edp}) gives $\rho \equiv 1$.  This suggests that we should seek for a coordinate transformation for adaptive mesh generation that is as close to the identity map as possible. Here,  we determine the coordinate transformation $\vx = \vx(\vxi)$  as the minimizer of the  $L^2$ cost functional \cite{Benamou:2000}
\begin{equation}\label{cost2} 
    C(\vx) = \int_{\Omega_c}\left|\vx(\vxi) - \vxi\right|^2 \mbox{d}\vxi
\end{equation} 
subject to the constraint (\ref{edp}).  From  \cite{Benamou:2000},  we find that the minimizer of the cost $C(\vx)$ in (\ref{cost2}) is the optimal solution of the $L^2$ optimal mass transfer  problem or $L^2$ Monge-Kantorovich problem (MKP).
%We would like to point out,  as noted in \cite{Benamou:200},
%the explicit dependence of $\vx$ and $\rho$ on the time variable  has been
%eliminated since the adaptive mesh or the optimal solution $\vx(\vxi)$ is
%computed at an arbitrary fixed point in time.  

In  \cite{Knott:1984, Brenier:1991}, it is shown that for bounded positive density function $\rho(\vx)$ and convex domains $\Omega_c$ and $\Omega$,   the solution of the $L^2$~MKP is unique and can be expressed as the gradient of some convex potential $\Psi$,  
\begin{equation}\label{AdpMesh}
  \vx(\vxi) = \nabla \Psi(\vxi) ,
\end{equation}
where  $\nabla$ is the gradient operator with respect to the computational variable $\mbox{\boldmath $\xi$}$.
Substituting  (\ref{AdpMesh}) into (\ref{edp}) we obtain the  Monge-Amp\`ere equation (MAE)
\begin{equation}\label{MA}
\left|\Omega_c\right|\rho\left(\nabla \Psi(\vxi)\right) \mbox{det}(D^2\Psi(\vxi)) = 1,
\end{equation}
where $\mbox{det} (D^2\Psi(\vxi))$ is the determinant of the Hessian of $\Psi$. 
   
As in  \cite{Sulman:2011a, Sulman:2011b},   we compute the solution of (\ref{MA}) as the steady-state solution of the parabolic Monge-Amp\`ere equation (PMA)
\begin{equation}\label{pmae}
\frac{\partial\Psi}{\partial \tau} = \log\left(\left|\Omega_c\right|\rho(\nabla\Psi) \mbox{det}\; D^2 \Psi\right),
\end{equation} 
with the initial and boundary conditions defined as 
\begin{equation}\label{incond}
 \Psi(\mbox{\boldmath $\xi$},0) = \Psi^0(\mbox{\boldmath $\xi$}) = \frac{1}{2}\mbox{\boldmath $\xi$}\cdot \mbox{\boldmath $\xi$}^T  
\end{equation} 
and 
\begin{equation}\label{nbd} 
\nabla \Psi \cdot \vn = \mbox{\boldmath $\xi$} \cdot \vn, \quad
\textrm{for}\;\; \mbox{\boldmath $\xi$} \in \partial \Omega_c,
\end{equation} 
where $\partial \Omega_c$ is the boundary of $\Omega_c$ and $\vn$ is the outward unit normal to $\partial \Omega_c$.  The boundary condition (\ref{nbd}) forces mesh points to stay on the boundary of the domain; they can only move along the  boundary. 

We would like to point out that if the solution of the physical model is time dependent, the initial condition (\ref{incond}) is employed only for computing  the initial adaptive mesh.  At the subsequent physical time levels,  the pseudo time integration of (\ref{pmae}) starts  at $\tau =0$ with the initial solution $\Psi^0$ taken as the steady-state solution obtained from the previous physical time level. 

The convergence of the solution of (\ref{pmae})--(\ref{nbd}) to the steady-state solution and the uniqueness
of the latter are shown  in \cite{Sulman:2011b}. 

Let $\Psi^\infty$ be the steady-state solution of (\ref{pmae}), (\ref{incond}), and (\ref{nbd}), then the adaptive mesh is determined  by taking the 
gradient of  $\Psi^\infty$, i.e.
\begin{equation}\label{adpx}
 \mbox{\boldmath $x$} = \nabla \Psi^\infty \, .
\end{equation}
To identify the steady-state solution of  (\ref{pmae}), we use the following stopping criterion:  
  \beq\label{dd_stop_criterion}
  \left\|\Psi^{n+1} - \Psi^n\right\|_2 = \left(\int_{\Omega_c} \left|\Psi^{n+1} - \Psi^n\right|^2 d \boldmath \xi\right)^{1/2} \le \mbox{TOL},
  \eeq
 where $\mbox{TOL}$ is the user specified tolerance.
 
We remark that the above procedure can be used to generate an adaptive mesh
for given analytical functions and steady-state and time-dependent problems.
To generate an adaptive mesh for a given function, starting from an initial mesh,
the monitor function is computed using the function value at the current mesh
and then the new mesh is generated by solving (\ref{pmae}) and (\ref{nbd}),
starting from the current mesh, until the steady state is reached.
To generate an adaptive mesh for a steady-state problem, the procedure is similar
except that in the current situation, the monitor function is calculated using
the computed solution on the current mesh and the physical model needs
to be re-solved on the new mesh for the new computed solution.
For a time-dependent problem, the monitor function is calculated based
on the computed solution and the mesh at the current time step and,
after the new mesh is obtained, the physical model is integrated
over one time step using the old and new meshes (see, for example, \cite{Huang2011}).

\section{The domain decomposition moving mesh method}\label{Sect3}

In this section, we describe a DD moving mesh method.  Motivated by the  Schwarz  methods,  we describe an overlapping domain decomposition technique for two dimensional domains.  The technique can be employed for domains in three spatial dimensions in an analogous way.  For simplicity, we consider the computational domain $\Omega_c = (0, \,1) \times (0, \, 1)$, and spilt  $\Omega_c$ into subdomains in one direction (i.e.,  either  in  $\xi$ or $\eta$ direction) or in two dimensions  (i.e., in both $\xi$ and $\eta$ directions). The  subdomains are obtained by the slab or block decompositions, respectively.

A slab decomposition in the $\xi$ direction  (the other slab or block decompositions are obtained in a similar manner) is obtained by  decomposing the $\xi$-interval $(0, 1)$ into $M$ subintervals in the $\xi$ direction, $(\alpha^i, \beta^i)$ for $i=1,\ldots, M$,  where  $\alpha^1 = 0$ and $\beta^M = 1$. The subdomains are required to overlap  in the following manner:
\begin{equation*}
 \alpha^i < \alpha^{i+1} < \beta^i  < \beta^{i+1}, \;\; \textrm{for}\;\;  i=1,\ldots, M.   
\end{equation*}
We then have $M$ subdomains $\Omega_i = (\alpha^i, \beta^i) \times (0, \, 1),$ for  $i=1,2, \ldots, M,$ in $\Omega_c$. 

\begin{figure}[H]
\centering
\includegraphics[width=12cm,height=6cm]{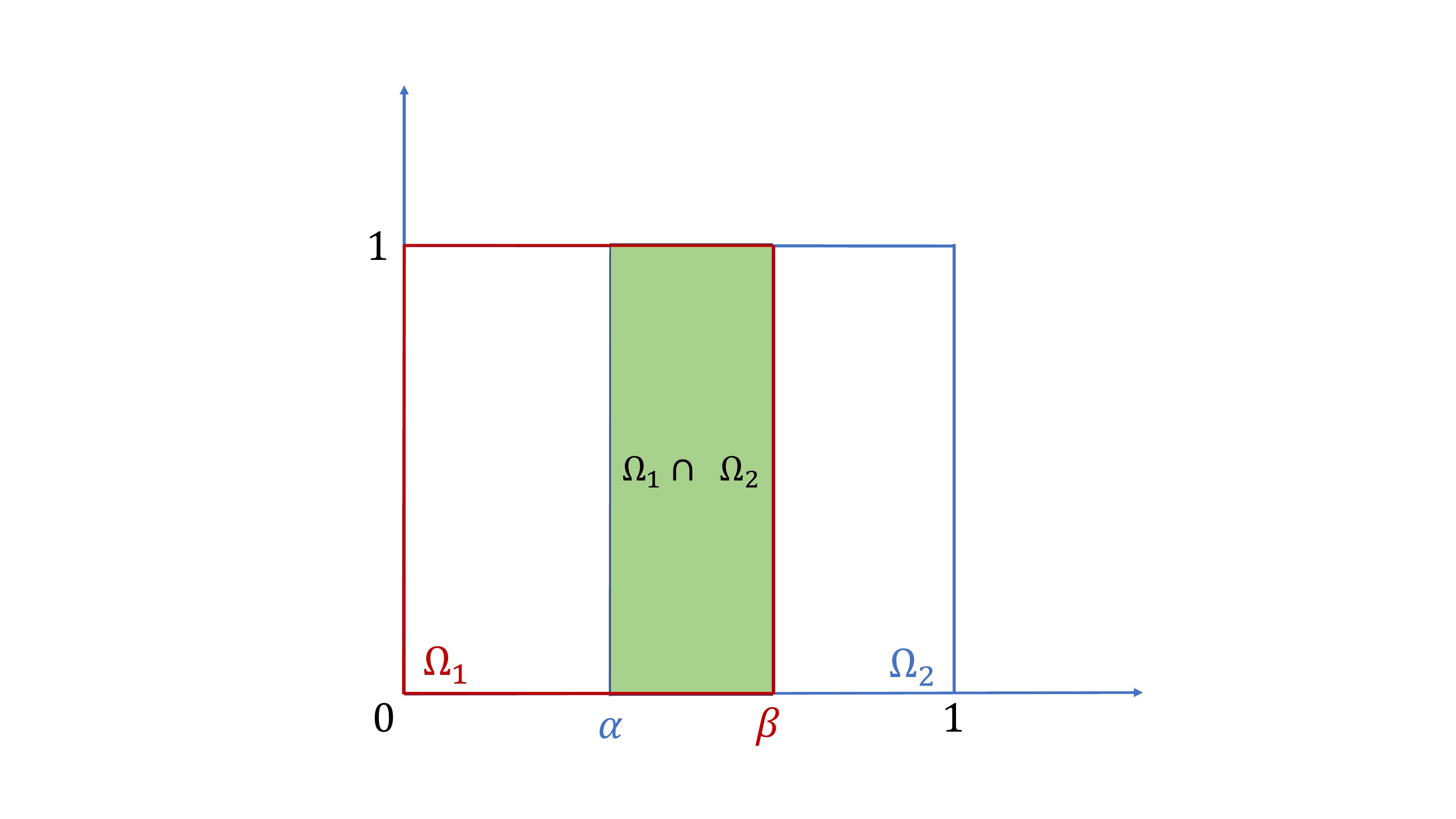}
\caption{\label{dd_pma_2subdomains} A domain decomposition of $\Omega_c$ in the $\xi$ direction into 2 subdomains $\Omega_1 = (0,\, \beta)\times(0,\,1)$ and $\Omega_2 = (\alpha,\, 1)\times(0,\,1),$ where $\alpha <\beta.$}
\end{figure}

For illustration purposes, we consider the case $M=2$, and divide $\Omega_c$ into two subdomains $\Omega_1$ and $\Omega_2$, as in Figure~\ref{dd_pma_2subdomains}.  Let $\Psi_i$ and $(x_i,y_i)$  be the solution of the Monge-Amp\`ere equation (\ref{MA}) and the corresponding coordinate transformation on the subdomain $\Omega_i$ for $i =1, 2.$  In this case, the DDPMA method computes $\Psi_1$ and $\Psi_2$ by solving  the two coupled initial value problems (IVPs): 
\begin{subequations}\label{pmae1}
 \begin{align}
 & \frac{\partial \Psi_1}{\partial \tau} = \log\left(\rho(\nabla\Psi_1) \mbox{det}\; D^2 \Psi_1\right), \quad \mbox{in}\quad \Omega_1,\label{pmae1a}\\[1em]
& \nabla \Psi_1(\xi,\eta, \tau) \cdot \vn = \mbox{\boldmath $\xi$} \cdot \vn, \quad \textrm{for}\;\; \mbox{\boldmath $\xi$} \in \partial \Omega_1 \cap \partial \Omega_c,\label{pmae1b}\\[1em]
& \Psi_1(\xi,\eta, \tau) = \Psi_2(\xi,\eta, \tau),\quad \mbox{on}\;\;   \partial \Omega_1 \cap  \overline{\Omega}_2,\label{pmae1c}\\[1em]
& \Psi_1(\mbox{\boldmath $\xi$},0) = \frac{1}{2}\mbox{\boldmath $\xi$}\cdot \mbox{\boldmath $\xi$}^T,  \quad \mbox{in}\quad \Omega_1,\label{pmae1d}
\end{align}
\end{subequations}

and
\begin{subequations}\label{pmae2}
 \begin{align}
 & \frac{\partial \Psi_2}{\partial \tau} = \log\left(\rho(\nabla\Psi_2) \mbox{det}\; D^2 \Psi_2\right), \quad \mbox{in}\quad \Omega_2,\label{pmae2a}\\[1em]
& \nabla \Psi_2(\xi,\eta,\tau) \cdot \vn = \mbox{\boldmath $\xi$} \cdot \vn, \quad \textrm{for}\;\; \mbox{\boldmath $\xi$} \in \partial \Omega_2 \cap \partial \Omega_c,\label{pmae2b}\\[1em]
& \Psi_2(\xi,\eta,\tau) = \Psi_1(\xi,\eta,\tau),\quad \mbox{on}\;\;   \partial \Omega_2 \cap  \overline{\Omega}_1,\label{pmae2c}\\[1em]
& \Psi_2(\mbox{\boldmath $\xi$},0) = \frac{1}{2}\mbox{\boldmath $\xi$}\cdot \mbox{\boldmath $\xi$}^T , \quad \mbox{in}\quad \Omega_2 .
\label{pmae2d}
\end{align}
\end{subequations}

We use standard centered finite differences for the spatial discretization of the parabolic Monge-Amp\`ere equations (\ref{pmae1}a) and (\ref{pmae2}a) to obtain the IVPs
\begin{subequations}\label{pmae1d}
 \begin{align}
 & \frac{d \Psi_1}{d \tau} = \log\left(\rho(\nabla_h\Psi_1) \mbox{det}\; D_h^2 \Psi_1\right), \quad \mbox{in}\quad \Omega_1,\label{pmae1a}\\[1em]
& \nabla_h \Psi_1(\xi,\eta, \tau) \cdot \vn = \mbox{\boldmath $\xi$} \cdot \vn, \quad \textrm{for}\;\; \mbox{\boldmath $\xi$} \in \partial \Omega_1 \cap \partial \Omega_c,\label{pmae1bd}\\[1em]
& \Psi_1(\xi,\eta, \tau) = \Psi_2(\xi,\eta,\tau),\quad \mbox{on}\;\;   \partial \Omega_1 \cap  \overline{\Omega}_2,\label{pmae1cd}\\[1em]
& \Psi_1(\mbox{\boldmath $\xi$},0) = \frac{1}{2}\mbox{\boldmath $\xi$}\cdot \mbox{\boldmath $\xi$}^T,  \quad \mbox{in}\quad \Omega_1,\label{pmae1dd}
\end{align}
\end{subequations}

and
\begin{subequations}\label{pmae2d}
 \begin{align}
 & \frac{d \Psi_2}{d\tau} = \log\left(\rho(\nabla_h\Psi_2) \mbox{det}\; D_h^2 \Psi_2\right), \quad \mbox{in}\quad \Omega_2,\label{pmae2a}\\[1em]
& \nabla_h \Psi_2(\xi,\eta,\tau) \cdot \vn = \mbox{\boldmath $\xi$} \cdot \vn, \quad \textrm{for}\;\; \mbox{\boldmath $\xi$} \in \partial \Omega_2 \cap \partial \Omega_c,\label{pmae2bd}\\[1em]
& \Psi_2(\xi,\eta,\tau) = \Psi_1(\xi,\eta,\tau),\quad \mbox{on}\;\;   \partial \Omega_2 \cap  \overline{\Omega}_1,\label{pmae2cd}\\[1 em]
& \Psi_2(\mbox{\boldmath $\xi$},0) = \frac{1}{2}\mbox{\boldmath $\xi$}\cdot \mbox{\boldmath $\xi$}^T,  \quad \mbox{in}\quad \Omega_2,\label{pmae2dd}
\end{align}
\end{subequations}
where $\nabla_h $ and $D_h^2 $ are the corresponding finite difference operators for the gradient and Hessian respectively. 

Notice that IVPs (\ref{pmae1d}) and (\ref{pmae2d}) are coupled through (\ref{pmae1cd}) and (\ref{pmae2cd}).
They  can be solved alternately or in parallel  for the steady state solutions to obtain $\Psi_1^\infty$ and $\Psi_2^\infty$, respectively. The coordinate transformations (and thus adaptive meshes) in $\Omega_1$ and $\Omega_2$ are  then determined by setting $\vx_1 = \nabla \Psi_1^\infty$ and  $\vx_2 = \nabla \Psi_2^\infty$, respectively. 

For the computation of the steady state solutions, we use an explicit variable time stepping time integrator for the IVPs (\ref{pmae1d}) and (\ref{pmae2d}) over small intervals of time $\Delta \tau$.   In our computation we use
Matlab function, ode113, a variable-step, variable-order Adams-Bashforth-Moulton PECE solver of orders 1 to 13
for this purpose.
Let  $\Psi_1^n$ and $\Psi_2^n$ be the solutions at time $\tau^n = n\Delta \tau$.  Let $(x_1^n,y_1^n) $ and $(x_2^n,y_2^n)$ be the corresponding adaptive meshes at time $\tau^n = n\Delta \tau$.  We summarize the steps for computing $\Psi_1^{n+1} $ and $\Psi_2^{n+1}$ together with the adaptive meshes $(x_1^{n+1},y_1^{n+1}) $ and $(x_2^{n+1},y_2^{n+1})$ alternately in  $\Omega_1$ and $\Omega_2$ in  Algorithm~\ref{alg1}.
\begin{algorithm}
\caption{(DDPMA method)}\label{alg1}
\begin{itemize}
\item[1.] Compute $\rho$ on $(x_1^n,y_1^n)$ and integrate  (\ref{pmae1d}) for one time step $\Delta\tau$ with the transmission condition 
\begin{equation}\label{transmission_condition1}
\Psi^{n+1}_1(\beta,\eta) = \Psi^n_2(\beta,\eta), \quad   \textrm{on}\;\; \partial \Omega_1 \cap  \overline{\Omega}_2
\end{equation}
 to obtain $\Psi_1^{n+1}$ on $\Omega_1$.
\item[2.] Compute the adaptive mesh  by setting $(x_1^{n+1},y_1^{n+1}) = \nabla \Psi_1^{n+1}$ .
\item[3.] Compute $\rho$ on $(x_2^n,y_2^n)$ and integrate  (\ref{pmae2d}) for one time step $\Delta\tau$ using the  transmission condition 
\begin{equation}\label{transmission_condition2}
\Psi^{n+1}_2(\alpha,\eta) = \Psi^{n+1}_1(\alpha,\eta),\quad \textrm{on}\;\; \partial \Omega_2 \cap  \overline{\Omega}_1
\end{equation}
 to obtain $\Psi_2^{n+1}$ on $\Omega_2$.
\item[4.] Compute the adaptive mesh  by setting $(x_2^{n+1},y_2^{n+1}) = \nabla \Psi_2^{n+1}$.
\item[5.] Compute $\mbox{res}_1 = \|\Psi_1^{n+1} - \Psi_1^n\|_2$ and $\mbox{res}_2 = \|\Psi_2^{n+1} - \Psi_2^n\|_2$.
\item[6.] Stop if  $\min\{\mbox{res}_1, \mbox{res}_2\} \le \mbox{TOL}$; Otherwise, set $n=n+1$ and go to 1.
\end{itemize}
\end{algorithm}
Notice that here we assume that $\Delta \tau$ is small enough so that the same mesh density function  $\rho(x^n,y^n)$
can be used for the time integration over the interval $(\tau^n,\tau^n+\Delta\tau)$. 
 
 In the special case of  $\Omega_c = (0,1)\times(0,1)$ and using slab decompositions,  the boundary conditions (\ref{pmae1b}) and (\ref{pmae2b})  can be expressed as
\beq\label{adapt_mesh_bc1}
x_1(0,\eta) = 0, \quad  0\le \eta \le 1,\quad  y_1(\xi,0) = 0,\;\;  y_1(\xi,1) = 1,\quad 0 \le \xi \le \beta,
\eeq
and
\beq\label{adapt_mesh_bc2}
x_2(1,\eta) = 1, \quad  0\le \eta \le 1,\quad  y_2(\xi,0) = 0,\;\;  y_2(\xi,1) = 1,\quad \alpha \le \xi \le 1.
\eeq
In this case, the transmission conditions  (\ref{transmission_condition1}) and (\ref{transmission_condition2})  (defined on the internal boundaries) take the form
\begin{equation}
\label{transmission_condition1x}
\Psi^{n+1}_1(\beta,\eta) = \Psi^n_2(\beta,\eta), \quad   0 \le \eta \le 1
\end{equation}
and 
\begin{equation}
\label{transmission_condition2x}
\Psi^{n+1}_2(\alpha,\eta) = \Psi^{n+1}_1(\alpha,\eta), \quad 0 \le \eta \le 1.
\end{equation}

\vs

A parallel DDPMA method is obtained by replacing  (\ref{transmission_condition2x}) by
\beq\label{transmission_condition3}
\Psi^{n+1}_2(\alpha,\eta) = \Psi^n_1(\alpha,\eta), \quad 0 \le \eta \le 1.
\end{equation}

Notice that this is a non-iterative domain decomposition algorithm - there is only one transfer of a subdomain solution information to its neighbouring subdomains per pseudo time step.  

\section{Numerical experiments}\label{Sect4}

In this section, we present several numerical experiments to demonstrate the performance of the DDPMA method described in Section~\ref{Sect3}.  We have used the DDPMA method with the alternating form of the transmission conditions (\ref{transmission_condition1x}) and (\ref{transmission_condition2x}) for serial computation and with the non-alternating form (\ref{transmission_condition3}) for parallel computation. We also give some numerical results on the convergence for the method. 

For the purpose of conducting these numerical experiments,  we choose the mesh density function $\rho(\vx)$ as the popular  arc-length function 
 \beq\label{monitfn}
  \rho(\vx) = \sqrt{1 + \vert \nabla_{\vx} u(\vx) \vert^2}, \quad \vx \in \Omega ,
  \eeq
where $u$ the solution of the physical model and $\nabla_{\vx}$ is the gradient operator with respect to $\vx$.

All the computations in Subsections~\ref{Sect4.1}  and~\ref{Sect4.2} have been done in  double precision Matlab on a mac computer with 2.3 GHz Intel Core i7 processor and 16 GB memory.

\subsection{A two-dimensional four subdomain decomposition}\label{Sect4.1} 

In this subsection, we present the results of the DDPMA method for computing  adaptive  meshes   in two spatial dimensions. The adaptive  mesh is computed with a 4-slab decomposition and $2 \times 2$ block decomposition, i.e.  4 subdomains in both cases. The subdomains are overlapping with an overlap of three grid points in both the slab and block decomposition cases. Notice that for the $2 \times 2$ block decomposition, the overlap  occurs in both the $\xi$-direction and $\eta$-direction. The alternating DDPMA method    described by Algorithem~\ref{alg1} is used for this example.

We consider two different examples of  the physical model. For these two examples, we assume the physical and computational domains are  $\Omega = \Omega_c = (0, 1)\times (0,1)$ and use a grid of size $65\times 65$ in the whole domain $\Omega_c$ to generate adaptive meshes in the physical domain $\Omega$.   The Matlab ODE solver, ode113, is used to integrate  the ODE systems (\ref{pmae1d}) and (\ref{pmae2d}) over each  time interval $(\tau^n, \tau^n+\Delta\tau)$ for $\Delta \tau=10^{-3}$. We have conducted a preliminary comparison with other ODE solvers and time steps and found that this choice of the ODE solver and time step gives better efficiency.

 In the first example, we  employ the DDPMA method to compute the adaptive mesh for the given function 
\begin{equation}\label{test_fn1}
u(x,y) = \fr{1}{ 1+ \exp({(x+y-1)}/{2\epsilon})}, \quad  \;\; (x,y) \in \Omega,
\end{equation}
which is an exact solution for the 2D Burgers' equation. We take $\epsilon=0.01$.

In the second example, we consider computing the  adaptive mesh for the function
        \beq\label{test_fn2}
        u(x,y) = 1 + \fr{9}{1 + 100r^2 \cos^2(\theta - 20r^2)}, \quad (x,y) \in \Omega,
        \eeq
        where   
        \[r = \sqrt{(x-0.7)^2 + (y - 0.5)^2},\quad \textrm{and}\quad \tan \theta = \fr{y-0.5}{x - 0.7}.\]
        
      Notice that the  Burgers' solution (\ref{test_fn1}) attains its maximum gradient along the line $x+y-1=0$, i.e. on the diagonal of the physical domain. Test function (\ref{test_fn2}) has its  maximum gradient along spiral shape that fills the whole physical domain. As a result,  the computed adaptive meshes are expected to concentrate along the line and the spiral shape, respectively.

\begin{figure}[h!]
\captionsetup{justification=justified,margin=.1cm}
\centering
\begin{subfigure}[]{0.3265\textwidth}
\includegraphics[width=1.75in,height=1.75in]{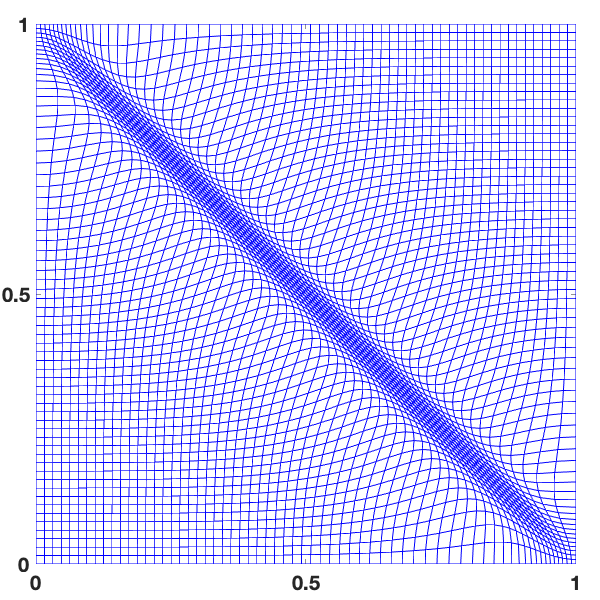}
\caption*{(a) single domain}
\end{subfigure}
\begin{subfigure}[]{0.3265\textwidth}
\includegraphics[width=1.75in,height=1.75in]{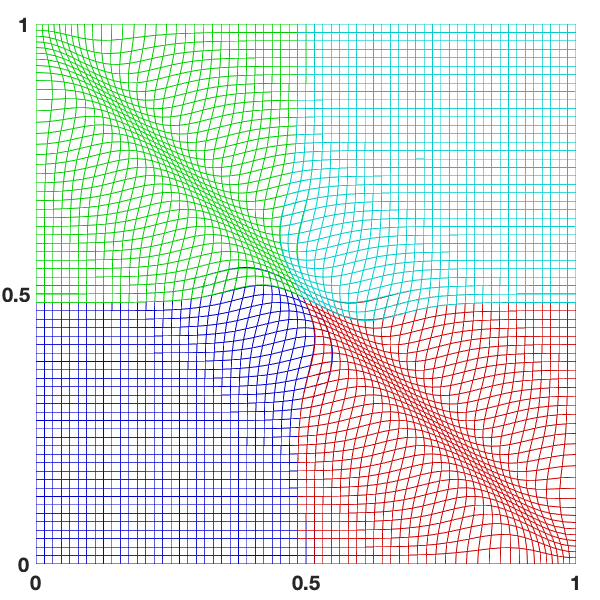}
\caption*{(b) 4-block}
\end{subfigure}
\begin{subfigure}[]{0.3265\textwidth}
\includegraphics[width=1.75in,height=1.75in]{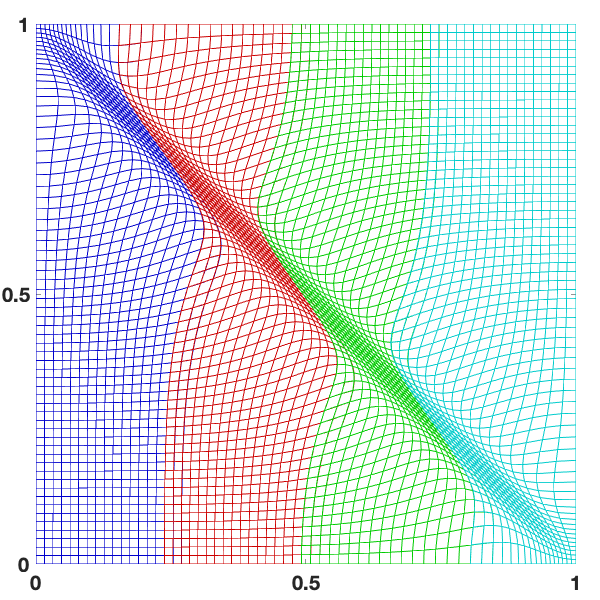}
\caption*{(c) 4-slab}
\end{subfigure}
\caption{\label{Burgers2dmesh} The  2D adaptive meshes computed for the test function (\ref{test_fn1}). Shown here are  adaptive meshes  obtained  by  (a) the PMA method employed on the whole domain,  (b) the DDPMA method using a 4-block decomposition, and (c) the DDPMA method using a 4-slab decomposition.}
\end{figure}

\begin{figure}[h!]
\captionsetup{justification=justified,margin=.1cm}
\centering
\begin{subfigure}[]{0.3265\textwidth}
\includegraphics[width=1.75in,height=1.75in]{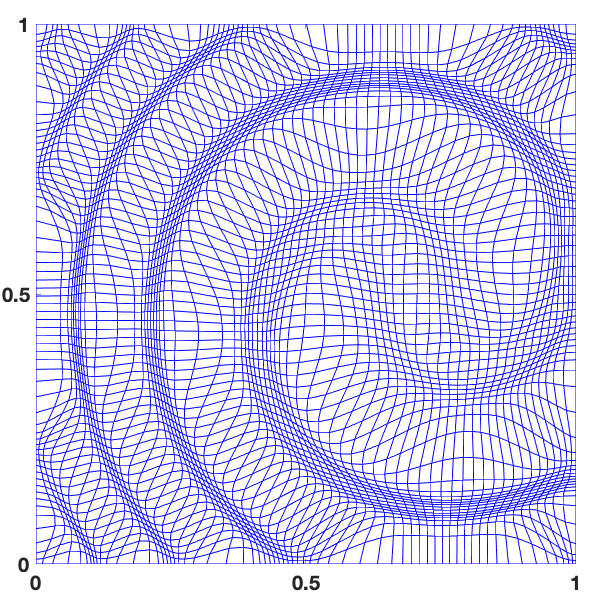}
\caption*{(a) Single domain}
\end{subfigure}
\begin{subfigure}[]{0.3265\textwidth}
\includegraphics[width=1.75in,height=1.75in]{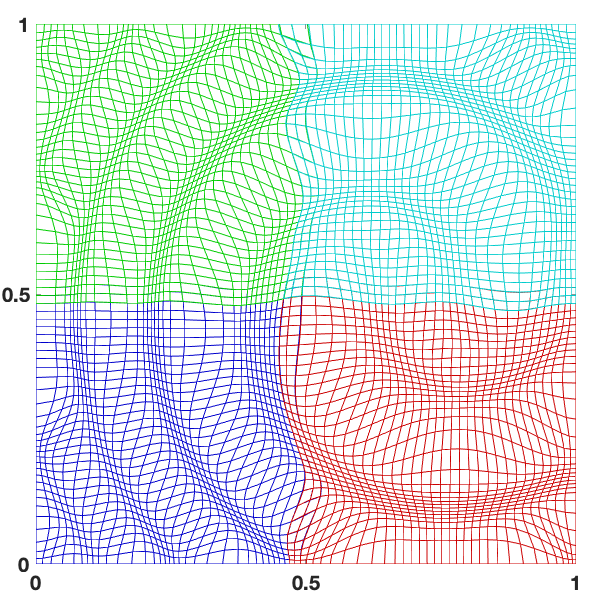}
\caption*{(b) 4-block}
\end{subfigure}
\begin{subfigure}[]{0.3265\textwidth}
\includegraphics[width=1.75in,height=1.75in]{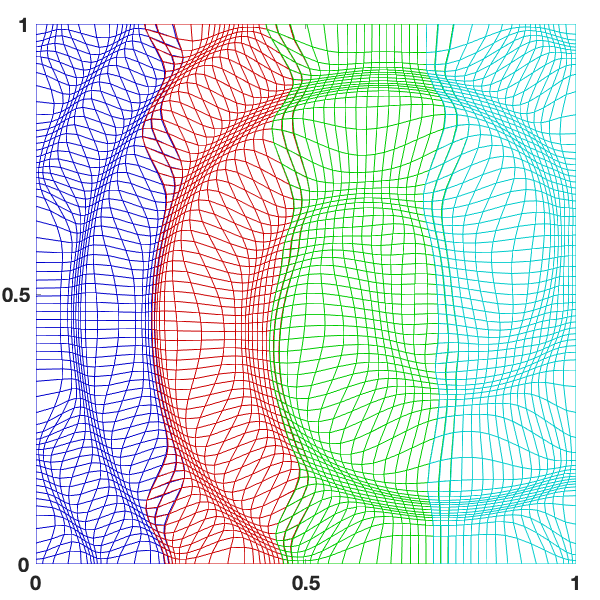}
\caption*{(c) 4-slab}
\end{subfigure}
\caption{\label{Spiral2dmesh}   The 2D adaptive meshes computed for the test function (\ref{test_fn2}). Shown here are  adaptive meshes  obtained by the  (a) the PMA method employed on the whole domain, (b) the DDPMA method using a 4-block decomposition, and (c)  DDPMA method using a  4-slab decomposition.}
\end{figure}

In Figure~\ref{Burgers2dmesh} (for function (\ref{test_fn1})) and Figure~\ref{Spiral2dmesh} (for function (\ref{test_fn2})) we present the adaptive meshes  computed  on a single domain and using a  4-block and 4-slab decomposition.  The figures show that there are no visible differences among the adaptive meshes computed by  the PMA method employed on the whole domain  and the DDPMA method using 4-slab and 4-block decompositions. 

 \subsection{A three-dimensional eight subdomain decomposition}\label{Sect4.2}
In this subsection, we illustrate the performance  of the DDPMA method for computing adaptive  meshes in three spatial dimensions. We  consider two different physical model solutions, and employ the DDPMA method to compute the adaptive mesh in  the physical domain. 
Similar to the 2D case, here  the  alternating DDPMA method with the transmission conditions (\ref{transmission_condition1}) and (\ref{transmission_condition2}) is implemented. We use a 8-block decomposition of the computational domain $\Omega_c$ for both examples. 

In the first example we assume that the physical solution is given as
\begin{equation}
\label{test_fn3d1}
u(x,y,z) = \tanh\left[100(x^2+y^2+z^2 -0.125)\right], \quad (x,y,z) \in (-1,1)^3.
\end{equation}
The solution $u(x,y,z)$ achieves its maximum gradient on the surface of a sphere centered at the origin with radius $r =\sqrt{.125}.$ Thus, the adaptive mesh is expected to be concentrated around the surface of the sphere. 

\begin{figure}[h!]
\captionsetup{justification=justified,margin=.3cm}
\centering
\includegraphics[width=6.75cm,height=6.75cm]{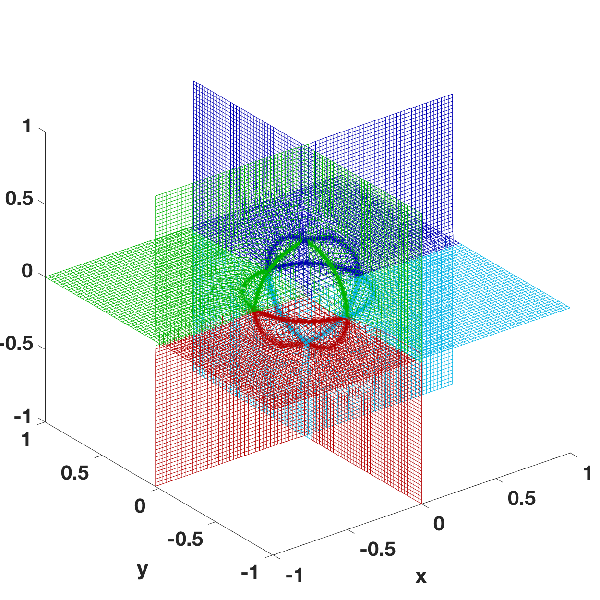}
\caption{\label{mesh3d1}  The 3D  adaptive mesh computed using the DDPMA method with 8-block decomposition for the test function (\ref{test_fn3d1}).}
\end{figure}

In Figure~\ref{mesh3d1} we show the adaptive meshes computed using the DDPMA method on three planes. The computed mesh appears to concentrate around a circle on each of the planes.  This illustrates that the computed mesh is concentrated around the surface of the sphere as expected. We  compare the adaptive mesh computed in each block with the mesh computed by applying the PMA method on  the whole domain, and we find excellent visual agreement between the two meshes which suggests convergence of the adaptive mesh obtained by the DDPMA method with four subdomains to the adaptive mesh obtained by the PMA method employed on a single domain. The CPU time for the alternating DDPMA method using 8 block is 274 seconds whereas it is  322 seconds for the PMA applied on the whole domain with a grid of size $81\times 81\times 81.$ 
  
\begin{figure}[h!]
\captionsetup{justification=justified,margin=.3cm}
\centering
       \begin{subfigure}[]{0.45\textwidth}
      \includegraphics[height=2in, width=2in]{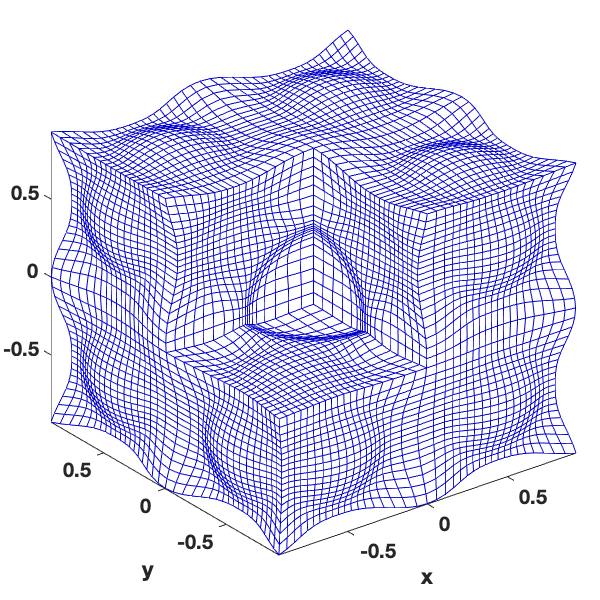}
\caption*{(a)~Single domain}
   \end{subfigure}
    \begin{subfigure}[]{0.45\textwidth}
 \centering    
      \includegraphics[height=2in, width=2in]{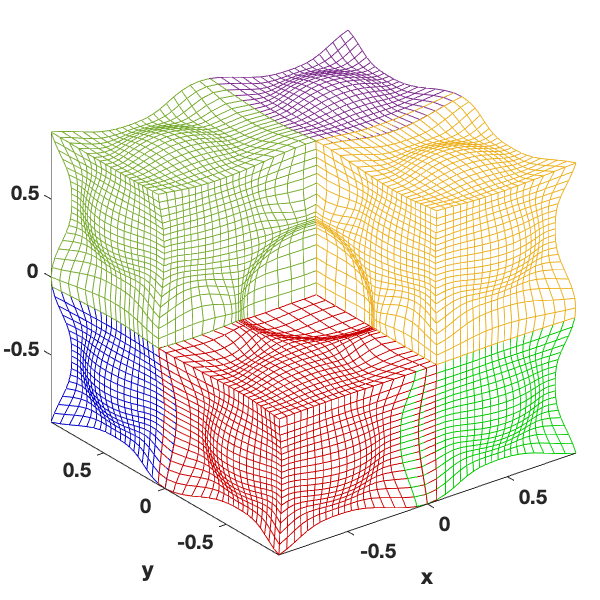} \qquad      
\caption*{(b)~8-block }
      \end{subfigure}           
   \caption{\label{mesh3d2} The 3D adaptive meshes computed  for the test function (\ref{test_fn3d2}). Shown are adaptive meshes computed using (a) the PMA method employed on the whole domain  and (b) the DDPMA method with 8-block decomposition.}
     \end{figure}

For the second example, we  employ the DDPMA method to compute the adaptive mesh for a test function that exhibits  sharp structures  that fill  the  whole physical domain.  To this end,  let $\ds{\Omega_c = \Omega= (-2,2)^3}$ and assume that the solution of the physical model is given as     
 \begin{equation}\label{test_fn3d2}
\ba{lcl}
\ds{u(x,y,z)} &=& \ds{\sum_{k=1}^9\tanh\left[50\left((x-x_0(k))^2 +  (y-y_0(k))^2  \right .\right .}\\[.5em]
 &&\qquad\qquad \ds{\left . \left .  +\, (z-z_0(k))^2 - 0.1875\right)\right],\quad (x,y,z)\in \Omega,}
\ea
\end{equation}
where
 \begin{align*}
 x_0 =  [0, 0.5, 0.5, -0.5, -0.5, 0.5, 0.5, -0.5, -0.5],\\
 y_0 =  [0, 0.5, -0.5, 0.5, -0.5, 0.5, -0.5, 0.5, -0.5],\\
 z_0 =  [0, 0.5, 0.5, 0.5, 0.5, -0.5, -0.5, -0.5, -0.5].
 \end{align*}
The test function (\ref{test_fn3d2}) has  its maximum gradients on the surfaces of nine spheres that are of radius  $r=\sqrt{0.1875}$ and centered at $(x_0(k), y_0(k), z_0(k))$, $k=1,2, \ldots, 9$ in the physical domain.  The adaptive mesh is expected to be clustered around  the surfaces of the nine spheres.  Figure~\ref{mesh3d2} presents adaptive meshes computed on a single domain and using an 8-block decomposition.  The concentration of the adaptive mesh can be seen along the surfaces of the nine spheres. The figure also shows a very good agreement of the adaptive mesh obtained by DDPMA method using 8-block decomposition and the adaptive mesh obtained by the PMA method employed on the whole domain.

\subsection{Efficiency of the DDPMA method}\label{Sect4.3}

We demonstrate the efficiency of the DDPMA method by examining the CPU time for the DDPMA method using four subdomains in 2D and eight subdomains in 3D.  We compare the results of the DDPMA method with the PMA method employed on the entire domain. To this end, we compute the adaptive meshes for different grid resolutions in both two and three spatial dimensions.  The DDPMA method is employed for both serial and parallel computations.   For the parallel algorithm,  the time integration of the ODE systems (\ref{pmae1d}) and (\ref{pmae2d}) from time level $n$ to $n+1$  on the subdomains is carried out in parallel using four processors for the case of 4 subdomains  and  8 processors for the case of 8 subdomains.

  The  computations in this subsection are conducted using double precision Matlab on a mac computer with 3.3 GHz 12-core Intel Xeon W  processor and 32 GB memory.

In Table~\ref{ddpma_efficiency2D} we present the CPU time required to compute the 2D adaptive meshes  for the test function (\ref{test_fn2}) using   the grid resolutions $65\times 65$, $129\times129$, $257\times 257$, $513\times513$ and $1025\times 1025$.   We can see  that for  two dimensional problems, the  DDPMA method becomes more efficient (in terms of the CPU time relative to the time required for the single domain solution)  as the number of the grid points increase.

Table~\ref{ddpma_efficiency3D} shows the CPU time to compute 3D adaptive meshes  for the test function (\ref{test_fn3d1}) using the grid resolutions  $65\times 65\times 65$, $81\times 81\times 81$, $101\times 101\times 101,$ and $121\times 121\times 121.$
 From these results it becomes clear that  the DDPMA method is more efficient than the PMA method employed on the entire domain.  The results presented here indicate that employing the DDPMA method in parallel improves the computational time significantly.  We would like to point out that one can employ the DDPMA method in parallel with more subdomains and  processors as needed to further speed up the computations.
 
\begin{table}[H]
\captionsetup{justification=justified,margin=.5cm}
  \centering
\begin{tabular}{|C{3cm}|C{2.5cm}|C{2.5cm}|C{2.5cm}|}
\hline
\multirow{2}{*}{grid size} & \multicolumn{3}{c|}{CPU time in seconds}   \\
\cline{2-4} &\multirow{2}{*}{ single domain} &  \multicolumn{2}{c|}{4-block subdomains}\\ 
\cline{3-4} &  & serial & parallel (4-core)\\
\hline
$65\times 65$ & 0.2 & 0.31 & 0.53 \\
\hline
 $ 129 \times 129$ & 0.38 &  0.51 & 0.66\\
 \hline
 $ 257 \times 257 $ & 0.8 &  1 & .76 \\
 \hline
 $ 513 \times 513$ & 3.53 & 2.7 & 1.6\\
\hline
 $1025\times 1025$ & 13.51 & 10.94 & 7.09\\
 \hline
\end{tabular}

\caption{\label{ddpma_efficiency2D} The 2D comparison of the CPU times  for the computation of  adaptive meshes using the DDPMA method with 4 subdomains and  the adaptive mesh computed using the PMA method on the entire domain. }
\end{table}

\begin{table}[H]
\captionsetup{justification=justified,margin=.5cm}
  \centering
\begin{tabular}{|C{3cm}|C{2.5cm}|C{2.5cm}|C{2.5cm}|}
\hline
\multirow{2}{*}{grid size} & \multicolumn{3}{c|}{CPU time in seconds}   \\
\cline{2-4} &\multirow{2}{*}{single domain} &  \multicolumn{2}{c|}{8-block subdomains}\\ 
\cline{3-4} &  & serial & parallel (8-core) \\
\hline
 $ 65 \times 65\times 65$ & 20 & 18 & 4\\
 \hline
 $ 81 \times 81\times 81$ & 30 & 29 & 8\\
 \hline
 $ 101 \times 101\times101$ & 62 & 58 & 16\\
 \hline
 $ 121 \times 121\times121$ & 104  & 94 & 37  \\
 \hline
\end{tabular}
 
\caption{\label{ddpma_efficiency3D} The 3D comparison of CPU times for the computation of adaptive meshes using  the DDPMA method with a 8-block decomposition and   the PMA method employed on the entire domain. The results are shown for computations conducted in   serial and parallel DDPMA with  eight processors. }
\end{table}

\subsection{Convergence of the DDPMA method}\label{Sect4.4}

In this subsection, we study the numerical convergence of the DDPMA method in two spatial dimensions.  We study the convergence of the solution obtained by the DDPMA method using the four-subdomain decomposition to the solution obtained by the PMA method on a single domain.

The analytical solution of the parabolic Monge-Amp\`ere equation (\ref{pmae}) is not available, therefore we use a solution computed by the PMA method on the single domain with a very fine grid resolution and for a very large number of pseudo time steps as a surrogate for the exact solution. 
Specifically, we consider the test function (\ref{test_fn2}) and employ the PMA method to compute the solution of the parabolic Monge-Amp\`ere equation (\ref{pmae}) with  a grid of size $1025\times 1025$ and 10000 pseudo time steps. This gives a solution $\Psi_{sd}$ that we assume is close enough to the exact solution.  Then, we employ the DDPMA method with 4 subdomains to  solve (\ref{pmae}) using the grid resolutions $33\times 33$, $65\times 65$, $129\times 129$ and $257\times 257$  and 1000 pseudo time steps to obtain a combined solution, $\Psi_{dd}$, on the union of the four subdomains.  The $L^p$ relative error is computed as
\beq\label{rel_errlp}
E_p = \fr{\left(\int_{\Omega_c} \left|\Psi_{sd} - \Psi_{dd}\right|^p d\xi d\eta\right)^{1/p}}{\left(\int_{\Omega_c} \left|\Psi_{sd}\right|^p d\xi d\eta\right)^{1/p}}
\eeq
for $p=1, 2$. The $L^\infty$ relative error is obtained as 
\beq\label{rel_errlinf}
E_\infty = \fr{\max \left|\Psi_{sd} - \Psi_{dd}\right|}{\max\left|\Psi_{sd}\right|} .
\eeq

\begin{figure}[h!]
\captionsetup{justification=justified,margin=.3cm}
\centering
\includegraphics[width=6.75cm,height=6.75cm]{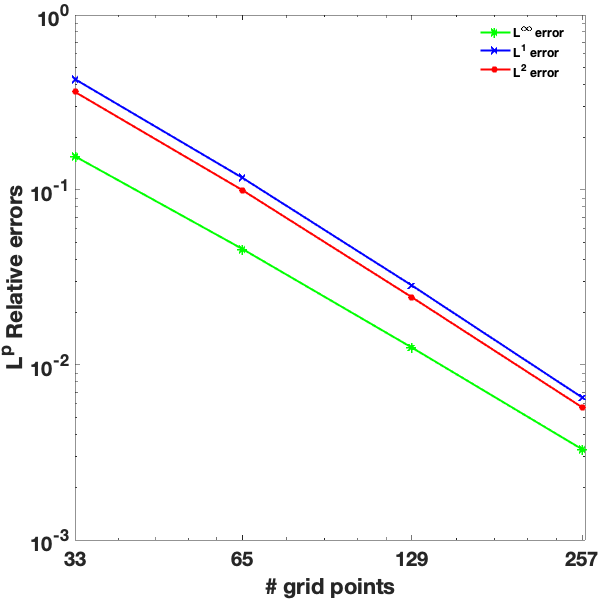}
\caption{\label{dd_convergence}  Convergence (in space) of the domain decomposition solution to the single domain solution. Shown are  the $L^p$ relative errors, for $p=1$, $2$, and $\infty$ (\ref{rel_errlp}) and (\ref{rel_errlinf}), for grid resolutions $N\times N$,  $N=33$, 65,  129, and 257.} 
\end{figure}

Figure~\ref{dd_convergence} presents the plots of the $L^p$ relative errors for $p=1$, 2, and $\infty$ in the logarithmic scale. The slopes of the $L^\infty$, $L^1$, and $L^2$ errors are 1.85, 2.01 and 1.99, respectively, which suggests that the convergence of the  DDPMA method is second order in space.  
 
\begin{figure}[h!]
\captionsetup{justification=justified,margin=.3cm}
\centering
\includegraphics[width=6.75cm,height=6.75cm]{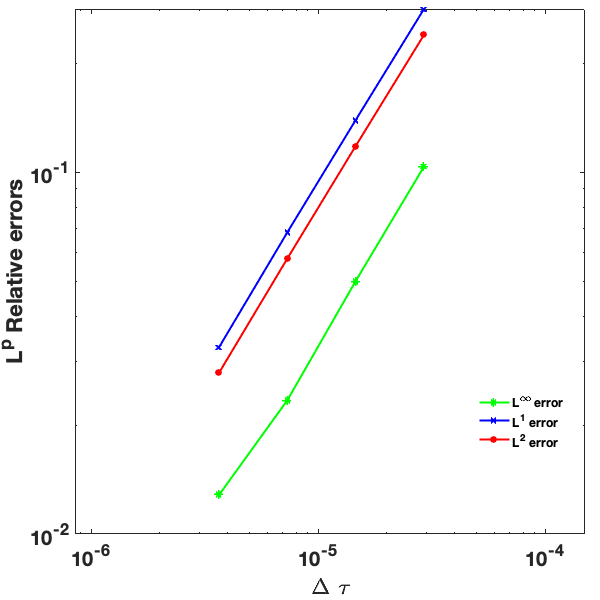}
\caption{\label{dd_convergencet}  Convergence rate of the domain decomposition solution to the single domain solution. Shown are  the $L^p$ relative errors, for $p=1$, 2, and $\infty$ (\ref{rel_errlp}) and (\ref{rel_errlinf}) plotted for the pseudo time steps $\Delta \tau = 0.015\Delta \xi \Delta \eta$, $0.03 \Delta \xi \Delta \eta$,  $0.06\Delta \xi \Delta \eta$, and $0.12 \Delta \xi \Delta \eta$.} 
\end{figure}

To examine the  convergence rate in the pseudo time variable. For the purpose of this test,  we  use the forward Euler method for the time integration of the ODE systems (\ref{pmae1d}) and (\ref{pmae2d}). To this end,  we fix the grid size to $65\times 65$ for whole domain and solve the parabolic Monge-Amp\`ere equation (\ref{pmae}) using small pseudo time step $\Delta \tau = 0.005\Delta \xi \Delta \eta$ and 5000 pseudo time steps to obtain the solution $\Psi_{sd}$. Then, we employ the DDPMA method with  four subdomains using same grid resolution but different pseudo time steps $\Delta \tau = 0.015\Delta \xi \Delta \eta$, $0.03 \Delta \xi \Delta \eta$,  $0.06\Delta \xi \Delta \eta$, and $0.12 \Delta \xi \Delta \eta$ to obtain $\Psi_{dd}$ for each pseudo time step size $\Delta\tau$.  We can then use the formulas (\ref{rel_errlp}) and (\ref{rel_errlinf}) to compute the $L^p$ relative errors for $p=1$, 2, and $\infty$.  In Figure~\ref{dd_convergencet}, we show the relative errors  in the logarithmic scale. For this test only, we use a forward Euler integrator in time, instead of the variable time stepping, variable order, ode113. We find that the slopes of  these relative errors are 1.01, 1.04 and 1.03, respectively. This illustrates that the convergence of the DDPMA method is first order in time.

\begin{figure}[h!]
\captionsetup{justification=justified,margin=.3cm}
\centering
       \begin{subfigure}[]{0.45\textwidth}
       \centering
        \includegraphics[height=2.5in, width=2.5in]{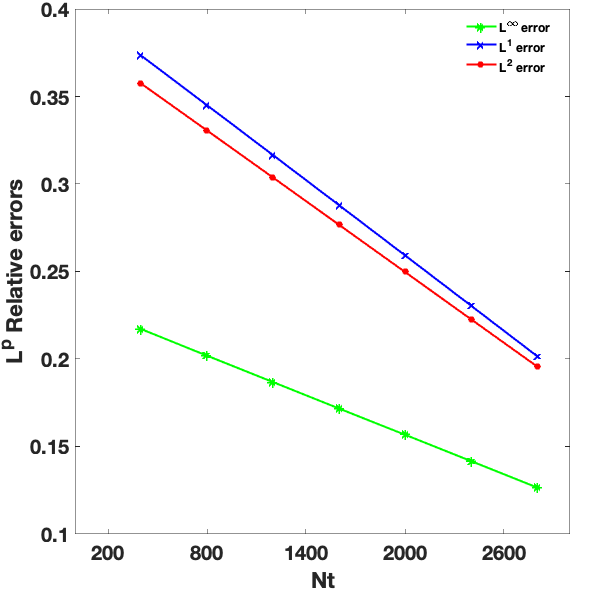}
   \end{subfigure} \qquad   
\caption{\label{dd_convergencetx}  The convergence of the domain decomposition solution to the single domain solution. Shown are  the $L^p$ relative errors, for $p=1$, 2, and $\infty$ (\ref{rel_errlp}) and (\ref{rel_errlinf}), as functions of the the number of pseudo time steps.} 
\end{figure}

We now study the convergence history of the domain decomposition solution as the number of time steps increases for a fixed mesh resolution and a fixed  $\Delta \tau$. We take the steady-state  single-domain solution of the same spatial resolution as the reference solution.  Figure~\ref{dd_convergencetx} shows the convergence history of the solution of the domain decomposition method to the single domain solution. It shows that the former converges to the latter as the number of pseudo time steps increases.

\begin{figure}[h!]
\captionsetup{justification=justified,margin=.3cm}
\centering
       \begin{subfigure}[]{0.45\textwidth}
       \centering
        \includegraphics[height=2.5in, width=2.5in]{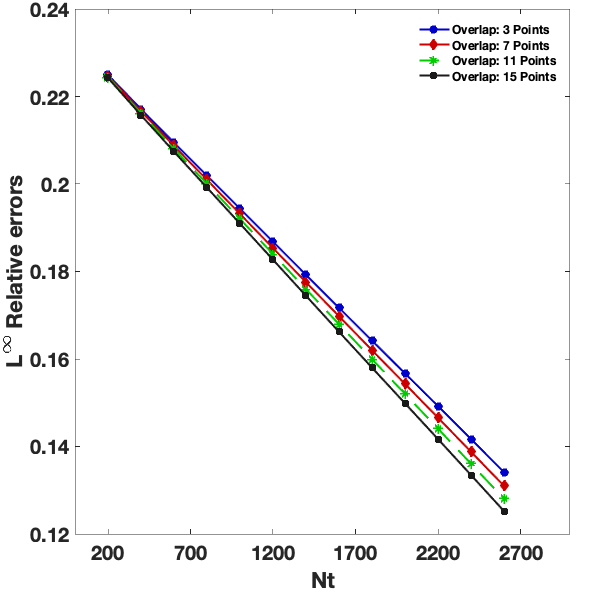}
   \end{subfigure} \qquad   
\caption{\label{dd_convergencetxx}  Convergence history of the domain decomposition solution to the single domain solution. Shown are  the $L^\infty$ relative error versus  the number of pseudo time steps for different number of overlap points.} 
\end{figure}

 To study the effect of the overlap of the subdomains on the convergence, in Figure~\ref{dd_convergencetxx} we plot the $L^\infty$ relative error (\ref{rel_errlinf}) versus 
the number of pseudo time steps for 5, 9, 11, and 15 overlap points. We notice that the error decreases as the number of the overlap points increases.

\begin{figure}[h!]
\captionsetup{justification=justified,margin=.3cm}
\centering
       \begin{subfigure}[]{0.45\textwidth}
       \centering
        \includegraphics[height=2.5in, width=2.5in]{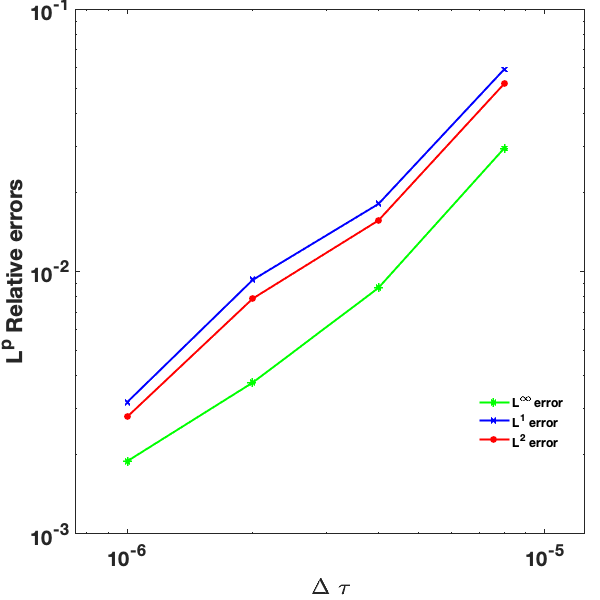}
   \end{subfigure} \qquad   
\caption{\label{dd_convergenceDt}  Convergence of the domain decomposition the steady state solution to the  single domain steady state solution. Shown are  the $L^\infty$ relative error versus  the time interval $\Delta \tau$.} 
\end{figure}

In what follows we study the convergence of the steady state solution obtained by the DDPMA algorithm to the steady state single 
domain solution $\Psi^\infty$ of the parabolic Monge-Amp\`ere equation (\ref{pmae}). We compute the DDPMA steady state solution using four subdomains for $\Delta \tau = 10^{-6}, 2\cdot 10^{-6}, 4\cdot 10^{-6}$ and $8\cdot 10^{-6}$.  A tolerance $\mbox{Tol} = 1e-6$ is used to detect the steady state solution; the time stepping is stopped when two successive solutions agree within the tolerance.  Figure~\ref{dd_convergenceDt} presents the plots of the $L^p$ relative errors versus  $\Delta \tau$. We notice here that as $\Delta \tau$ gets smaller the relative errors decrease which illustrates the convergence of the DDPMA steady state solution to the steady state solution obtained on a single domain.

\subsection{Quality measures of the DDPMA adaptive meshes}\label{Sect4.5}

In this subsection, we compute the quality measure of the adaptive meshes computed using the DDPMA method.  In $2D$, the adaptive mesh on the physical domain $\Omega$ is formed by combining the adaptive meshes computed on each  of  the four subdomains $\Omega_i, \;\;  i = 1, 2, 3, 4. $  On each subdomain $\Omega_i$,  the adaptive mesh  is obtained as an image of a coordinate transformation computed using the PMA method.  We use a mesh quality measure $E_{\mbox{adp}}$ as described in \cite{Sulman:2011a,Huang2005} which is given by  
\begin{equation}
\label{quality_measure}
E_{\mbox{adp}}(\vx) = \frac{\rho(\nabla \Psi^{\infty})\mbox{J}}{|\Omega_c|}, \quad \forall \vx \in \Omega ,
\end{equation}
 where $\Psi^{\infty}$   is the DD solution obtained by combining  the  steady state solutions $\Psi^{\infty}_i$ of the parabolic Monge-Amp\`ere equation on the subdomains $\Omega_i$, $|\Omega_c|$ is the area (in 2D) or volume (in 3D) of the computational domain, and $\mbox{J}$ is the determinant of the Jacobian matrix of the coordinate transformation. Note that $E_{\mbox{adp}}$ is defined pointwise for each grid node in the domain $\Omega$. We compute the maximum and $L^2$ norms of the  mesh quality measure, $E_{max} = ||E_{\mbox{adp}}||_{\infty}$ and $E_2=||E_{\mbox{adp}}||_2$, respectively.

The mesh  quality measure (\ref{quality_measure})  is computed for the adaptive mesh obtained by the DDPMA method using 4 slabs and   a $2 \times 2$ block  decomposition. Here, the mesh quality  measure $E_{\mbox{adp}}$ is  computed  for two different grid resolutions, namely using $21 \times 21$ and $41 \times 41$ mesh points. 

The results from  Table \ref{quality_measure_table} show  that the  adaptive meshes computed using the DDPMA method  and PMA method on a single domain have very similar mesh quality measures. This indicates excellent agreement between the adaptive mesh obtained using  the DDPMA and the adaptive mesh  obtained using the PMA method employed on the entire domain. Moreover, the fact that the values in the table are close to one indicates that the meshes satisfy the equidistribution principle (\ref{edp}) closely.

\begin{table}[H]
\centering
\begin{tabular}{ |C{5cm}|C{1.5cm}|C{1.5cm}|C{1.5cm}|C{1.5cm}| }
\hline
\multirow{2}{*}{\diagbox{\begin{scriptsize} Decompositions \end{scriptsize}}{\begin{scriptsize} Qual. measure \end{scriptsize} }} & \multicolumn{2}{c|}{$21 \times 21$} & \multicolumn{2}{c|}{$41 \times 41$} \\
\cline{2-5}
& $E_2 $ & $E_{max} $ & $E_2 $ & $E_{max} $  \\
\hline
4 Slabs DD & $1.0198$ & $1.1546$ & $1.0053$ & $1.0997$   \\
\hline
$2\times 2$ Block DD & $1.0194$ & $1.1590$ & $1.0052$ & $1.0994$  \\
\hline
Whole Domain & $1.0194$ & $1.1609$ & $1.0042$ & $1.0980$  \\
\hline
\end{tabular}
\caption{\label{quality_measure_table} Mesh quality measure of the adaptive meshes generated by the DDPMA and PMA methods.}
\end{table}

\section{Conclusions} \label{Sect5}

We have developed a non-iterative  overlapping domain decomposition approach for fast and efficient computation of   adaptive moving meshes in multi-dimensions. The computational domain is split into subdomains and  the parabolic Monge-Amp\`ere method is employed to compute the adaptive mesh on each subdomain.  The numerical experiments show that the DDPMA  method is more efficient than the PMA method applied on the whole 
domain. This result is significant especially when solving physical problems on large 2D domains and/or in three spatial dimensions.
 The computations involved here are performed on both a single processor (for the serial computations) and 4-processors (for the parallel computations). The number of the processors utilized can be increased by increasing the number of subdomains.  The results indicate that the parallel computations can be  implemented efficiently with the DDPMA method.  
We have also studied the convergence of the adaptive mesh computed using the DDPMA method to the adaptive mesh computed by the PMA method applied to the whole domain.

\section*{References}
  \bibliographystyle{elsarticle-num} 
    \bibliography{ddpma}

\begin{thebibliography}{10}
\expandafter\ifx\csname url\endcsname\relax
  \def\url#1{\texttt{#1}}\fi
\expandafter\ifx\csname urlprefix\endcsname\relax\def\urlprefix{URL }\fi
\expandafter\ifx\csname href\endcsname\relax
  \def\href#1#2{#2} \def\path#1{#1}\fi

\bibitem{Dorfi:1987}
E.~Dorfi, L.~Drury, Simple adaptive grids for 1-{D} initial value problems, J.
  Comput. Phys. 69~(1) (1987) 175--195.

\bibitem{Thompson:1985}
J.~F. Thompson, Z.~U.~A. Warsi, C.~W. Mastin, Numerical {G}rid {G}eneration:
  {F}oundations and {A}pplications, North-Holland Publishing Co., New York,
  1985.

\bibitem{Anderson:1987}
D.~A. Anderson, Equidistribution schemes, {P}oisson generators, and adaptive
  grids, Appl. Math. Comp. 24~(3) (1987) 211--227.

\bibitem{Huang2011}
W.~Huang, R.~D. Russell, Adaptive Moving Mesh Methods, Springer, New York, USA,
  2011.

\bibitem{Miller:1981}
K.~Miller, Moving finite elements {II}, SIAM J. Numer. Anal. 18~(6) (1981)
  1033--1057.

\bibitem{Millera:1981}
K.~Miller, R.~N. Miller, Moving finite elements {I}, SIAM J. Numer. Anal.
  18~(6) (1981) 1019--1032.

\bibitem{Gelinas:1981}
R.~Gelinas, S.~Doss, K.~Miller, The moving finite element method: Applications
  to general partial differential equations with multiple large gradients, J.
  Comput. Phys. 40~(1) (1981) 202--249.

\bibitem{Adjerid:1986}
S.~Adjerid, J.~E. Flaherty, A moving finite element method with error
  estimation and refinement for one-dimensional time dependent partial
  differential equations, SIAM J. Numer. Anal. 23~(4) (1986) 778--796.

\bibitem{Furrzeland:1990}
R.~Furzeland, J.~Verwer, P.~Zegeling, A numerical study of three moving-grid
  methods for one-dimensional partial differential equations which are based on
  the method of lines, J. Comput. Phys. 89~(2) (1990) 349--388.

\bibitem{Hawken:1991}
D.~Hawken, J.~Gottlieb, J.~Hansen, Review of some adaptive node-movement
  techniques in finite-element and finite-difference solutions of partial
  differential equations, J. Comput. Phys. 95~(2) (1991) 254--302.

\bibitem{Adjerid:1992}
S.~Adjerid, J.~E. Flaherty, P.~K. Moore, Y.~J. Wang, High-order adaptive
  methods for parabolic systems, Phys. D 60~(1-4) (1992) 94--111.

\bibitem{Huang:1994}
W.~Huang, Y.~Ren, R.~D. Russell, Moving mesh partial differential equations
  ({MMPDES}) based on the equidistribution principle, SIAM J. Numer. Anal.
  31~(3) (1994) 709--730.

\bibitem{Huang:1999}
W.~Huang, R.~D. Russell, Moving mesh strategy based on a gradient flow equation
  for two-dimensional problems, SIAM J. Sci. Comput. 20~(3) (1999) 998--1015.

\bibitem{Huang:2003}
W.~Huang, W.~Sun, Variational mesh adaptation {II}: error estimates and monitor
  functions, J. Comput. Phys. 184~(2) (2003) 619--648.

\bibitem{Cao:2002}
W.~Cao, W.~Huang, R.~D. Russell, A moving mesh method based on the geometric
  conservation law, SIAM J. Sci. Comput. 24~(1) (2002) 118--142.

\bibitem{Backett:2002}
G.~Beckett, J.~A. Mackenzie, A.~Ramage, D.~M. Sloan, Computational solution of
  two-dimensional unsteady {PDE}s using moving mesh methods, J. Comput. Phys.
  182~(2) (2002) 478--495.

\bibitem{Boor:1973}
C.~de~Boor, Good approximation by splines with variable knots. {II}, in:
  Conference on the Numerical Solution of Differential Equations (Univ. Dundee,
  Dundee, 1973), Springer, Berlin, 1974, pp. 12--20. Lecture Notes in Math.,
  Vol. 363.

\bibitem{Huang:1994b}
W.~Huang, Y.~Ren, R.~D. Russell, Moving mesh methods based on moving mesh
  partial differential equations, J. Comput. Phys. 113~(2) (1994) 279--290.

\bibitem{Tan05}
T.~Tang, Moving mesh methods for computational fluid dynamics flow and
  transport, in: Recent Advances in Adaptive Computation (Hangzhou, 2004), Vol.
  383 of AMS Contemporary Mathematics, Amer. Math. Soc., Providence, RI, 2005,
  pp. 141--173.

\bibitem{BHR09}
C.~J. Budd, W.~Huang, R.~D. Russell, Adaptivity with moving grids, Acta
  Numerica 18 (2009) 111--241.

\bibitem{Budd:2009}
C.~J. Budd, J.~F. Williams, Moving mesh generation using the parabolic
  {M}onge-{A}mp\`ere equation., SIAM J. Sci. Comput. 31~(5) (2009) 3438--3465.

\bibitem{Sulman:2011a}
M.~Sulman, J.~F. Williams, R.~D. Russell, Optimal mass transport for higher
  dimensional adaptive grid generation, J. Comput. Phys. 230~(9) (2011)
  3302--3330.

\bibitem{Chacon:2011}
L.~Chac{\'o}n, G.~Delzanno, J.~Finn, Robust, multidimensional mesh-motion based
  on {M}onge-{K}antorovich equidistribution, J. Comput. Phys. 230~(1) (2011)
  87--103.

\bibitem{Monge:1781}
G.~Monge, M\'emoire sur la th\'eorie des d\'eblais at des remblais, in:
  Histoire de l'Acad\'emie Royale des Sciences de Paris, 1781, pp. 666--704.

\bibitem{Kantorovich:1948}
L.~V. Kantorovich, On a problem of {M}onge, Uspehki Mat. Nauk 3 (1948)
  225--226.

\bibitem{Benamou:2000}
J.-D. Benamou, Y.~Brenier, A computational fluid mechanics solution to the
  {M}onge-{K}antorovich mass transfer problem, Numer. Math. 84~(3) (2000)
  375--393.

\bibitem{Sulman:2011b}
M.~Sulman, J.~Williams, R.~D. Russell, An efficient approach for the numerical
  solution of the {M}onge-{A}mp\`ere equation, Appl. Numer. Math. 61~(3) (2011)
  298--307.

\bibitem{Haynes:2012}
M.~Gander, R.~Haynes, Domain decomposition approaches for mesh generation via
  the equidistribution principle, SIAM J. Numer. Anal. 50 (2012) 2111--2135.

\bibitem{Haynes:2017}
R.~D. Haynes, F.~Kwok, Discrete analysis of domain decomposition approaches for
  mesh generation via the equidistribution principle, Math. Comp. 86~(303)
  (2017) 233--273.

\bibitem{Haynes:2014}
R.~Haynes, A.~Howse, Generating equidistributed meshes in {2D} via domain
  decomposition, in: Domain Decomposition Methods in Science and Engineering
  XXI, Vol.~98 of Lecture Notes in Computational Science and Engineering,
  Springer, 2014, pp. 167--178.

\bibitem{Bihlo:2014}
A.~Bihlo, R.~D. Haynes, Parallel stochastic methods for {PDE} based grid
  generation, Comput. Math. Appl. 68~(8) (2014) 804--820.

\bibitem{Haynes:2018}
R.~D. Haynes, Domain decomposition approaches for {PDE} based mesh generation,
  in: Domain decomposition methods in science and engineering {XXIV}, Vol. 125
  of Lect. Notes Comput. Sci. Eng., Springer, Cham, 2018, pp. 73--86.

\bibitem{Schwarz:1869}
H.~A. Schwarz, Uber einige abbildungsaufgaben, Ges. Math. Abh. 11 (1869)
  65--83.

\bibitem{Lions:1989}
P.~Lions, On the {S}chwarz alternating method {I}, in: T.~F. Chan, et~al.
  (Eds.), Domain Decomposition Methods, SIAM, Philadelphia, 1989, Ch.~10.

\bibitem{Babuska:1958}
I.~Babuska, On the {S}chwarz algorithm in the theory of differential equations
  of mathematical physics, Tchecosl. Math J. 8 (1958) 328--342.

\bibitem{Michlin:1951}
S.~G. Michlin, On the {S}chwarz algorithm, Dokl. Acad. N. USSR. 77 (1951)
  569--571.

\bibitem{White:1987}
J.~K. White, A.~Sangiovanni-Vincentelli, Waveform relaxation, in: Relaxation
  Techniques for the Simulation of VLSI Circuits, Springer, Boston, MA, 1987,
  pp. 79--100.

\bibitem{Jeltsch:1995}
R.~Jeltsch, B.~Pohl, Waveform relaxation with overlapping splittings, SIAM J.
  Sci. Comput., 16~(1) (1995) 40--49.

\bibitem{Gander:1997}
M.~J. Gander, Overlapping {S}chwarz waveform relaxation for parabolic problems,
  in: In Proceedings of Algoritmy'97, 1997, pp. 425--431.

\bibitem{Gander:2002}
M.~J. Gander, H.~Zhao, Overlapping {S}chwarz waveform relaxation for the heat
  equation in n dimensions, BIT Numer. Math. 42~(4) (2002) 779--795.

\bibitem{Gander:2005}
S.~Vandewalle, M.~J. Gander, Optimized overlapping {S}chwarz methods for
  parabolic {PDE}s with time-delay, in: T.~J. Barth, M.~Griebel, D.~E. Keyes,
  R.~M. Nieminen, D.~Roose, T.~Schlick, R.~Kornhuber, R.~Hoppe, J.~P{\'e}riaux,
  O.~Pironneau, O.~Widlund, J.~Xu (Eds.), Domain Decomposition Methods in
  Science and Engineering, Springer Berlin Heidelberg, Berlin, Heidelberg,
  2005, pp. 291--298.

\bibitem{Vabishchevich:2008}
P.~Vabishchevich, Domain decomposition methods with overlapping subdomains for
  the time-dependent problems of mathematical physics, Comput. Methods Appl.
  Math. 8~(4) (2008) 393--405.

\bibitem{Cai:1994}
X.-C. Cai, Multiplicative {S}chwarz methods for parabolic problems, SIAM J.Sci.
  Comput. 15~(3) (1994) 587--603.

\bibitem{Zheng:2008}
Z.~Zheng, B.~Simeon, L.~Petzold, A stabilized explicit {L}agrange multiplier
  based domain decomposition method for parabolic problems, J. Comput. Phys.
  227~(10) (2008) 5272--5285.

\bibitem{Qin:2008}
L.~Qin, X.~Xu, Optimized {S}chwarz methods with {R}obin transmission conditions
  for parabolic problems, SIAM J. Sci. Comput. 31~(1) (2008) 608--623.

\bibitem{Mota:2017}
A.~Mota, I.~Tezaur, C.~Alleman, The {S}chwarz alternating method in solid
  mechanics, Comput. Meth. Appl. Mech. Eng. 319~(1) (2017) 19--51.

\bibitem{Gnatyuk:2015}
M.~A. Gnatyuk, V.~M. Morozov, On the {S}chwarz alternating method for solving
  electromagnetic problems, in: 2015 XXth IEEE International Seminar/Workshop
  on Direct and Inverse Problems of Electromagnetic and Acoustic Wave Theory
  (DIPED), Lviv, Ukraine, 2015, pp. 132--135.

\bibitem{Dai:2016}
Z.~Dai, Q.~Du, B.~Liu, Schwarz alternating methods for anisotropic problems
  with prolate spheroid boundaries, SpringerPlus 5~(1423) (2016) 1.

\bibitem{Dawsonfinitedifferencedomain1991}
C.~N. Dawson, Q.~Du, T.~F. Dupont, A finite difference domain decomposition
  algorithm for numerical solution of the heat equation, Math. Comput. 57~(195)
  (1991) 63--71.

\bibitem{DawsonExplicitImplicitConservative1994}
C.~Dawson, T.~Dupont, Explicit/implicit, conservative domain decomposition
  procedures for parabolic problems based on block-centered finite differences,
  SIAM J. Numer. Anal. 31~(4) (1994) 1045--1061.

\bibitem{ZhangStablegloballynoniterative2000}
Y.~Zhang, Stable, globally non-iterative, non-overlapping domain decomposition
  methods for the efficient solution of parabolic evolutionary systems, Ph.D.
  thesis, Louisiana State University (2000).

\bibitem{GuangweiUNCONDITIONALSTABILITYPARALLEL2007}
G.~Yuan, Z.~Sheng, X.~Hang, The unconditional stability of parallel difference
  schemes with second order convergence for nonlinear parabolic system, J.
  Partial Diff. Eq. 20 (2007) 45--64.

\bibitem{YangNoniterativeparallelSchwarz2017}
D.~Yang, Non-iterative parallel {{Schwarz}} algorithms based on overlapping
  domain decomposition for parabolic partial differential equations, Math.
  Comput. 86~(308) (2017) 2687--2718.

\bibitem{Xuenewparallelalgorithm2018}
G.~Xue, H.~Feng, A new parallel algorithm for solving parabolic equations, Adv.
  Diff. Eq. 2018~(174) (2018) 1--16.

\bibitem{Knott:1984}
M.~Knott, C.~S. Smith, On the optimal mapping of distributions, J. Optim.
  Theory Appl. 43~(1) (1984) 39--49.

\bibitem{Brenier:1991}
Y.~Brenier, Polar factorization and monotone rearrangement of vector-valued
  functions, Comm. Pure Appl. Math. 44~(4) (1991) 375--417.

\bibitem{Huang2005}
W.~Huang, Measuring mesh qualities and application to variational mesh
  adaptation, SIAM J. Sci. Comput. 26 (2005) 1643--1666.

\end{thebibliography}

\end{document}